\documentclass[12pt]{article}
\input epsf
\usepackage{amsmath}
\usepackage{amssymb}
\usepackage{theorem}

\numberwithin{equation}{section}

\newtheorem{thm}{Theorem}[section]
\newtheorem{lemma}[thm]{Lemma}
\newtheorem{prop}[thm]{Proposition}
\newtheorem{cor}[thm]{Corollary}
{\theorembodyfont{\rmfamily}
\newtheorem{defn}[thm]{Definition}
\newtheorem{rmk}[thm]{Remark}}

\newcommand{\qed}{\hfill \mbox{\raggedright \rule{.07in}{.1in}}}
 
\newenvironment{proof}{\vspace{1ex}\noindent{\bf
Proof}\hspace{0.5em}}{\hfill\qed\vspace{1ex}}

\newcommand{\R}{{\mathbb R}}
\newcommand{\T}{{\mathbb T}}
\newcommand{\Z}{{\mathbb Z}}
\newcommand{\N}{{\mathbb N}}

\renewcommand{\hat}{\widehat}
\renewcommand{\tilde}{\widetilde}
\renewcommand{\bar}{\overline}

 \newcommand{\diam}{\operatorname{diam}}

\newcommand{\SMALL}{\textstyle}

\title{Almost Sure Invariance Principle for Nonuniformly Hyperbolic Systems}
\author{
Ian Melbourne 
\\ Department of Maths and Stats
\\ University of Surrey
\\ Guildford GU2 7XH, UK
\and
Matthew Nicol  
\\ Department of Maths 
\\ University of Houston
\\ Houston TX 77204-3008, USA
}

\date{23 September, 2004.   Revised 8 February, 2005.}

\begin{document}

\maketitle

\begin{abstract}
We prove an almost sure invariance principle that is valid for general
classes of nonuniformly expanding and nonuniformly hyperbolic dynamical 
systems.    Discrete time systems and flows are covered by this result.
In particular, the result applies to the planar periodic Lorentz 
flow with finite horizon.

Statistical limit laws such as the central limit theorem,
the law of the iterated logarithm, and their functional versions,
are immediate consequences.

\end{abstract}

\section{Introduction}

Statistical properties of uniformly expanding maps and uniformly
hyperbolic (Axiom~A) diffeomorphisms are by now classical.
H\"older observations satisfy exponential decay
of correlations and the central limit theorem (CLT), see for example 
Bowen~\cite{Bowen75}, Ratner~\cite{Ratner73}, Ruelle~\cite{Ruelle78},
Parry and Pollicott~\cite{ParryPoll90}.
Furthermore,
Denker and Philipp~\cite{DenkerPhilipp84} proved an almost sure invariance
principle (ASIP) for H\"older observations.    Immediate consequences of
the ASIP are the CLT, the law of the iterated logarithm (LIL),
and their functional versions, see~\cite{PhilippStout75}.

Many proofs of the CLT for dynamical systems use directly
the martingale approximation method of Gordin~\cite{Gordin69}, 
see~\cite{Keller80,Liverani96,MT02}.   The ASIP
can often also be obtained in this way
see~\cite{ConzeBorgne01,FMT03,MT02,Walkden} and indeed this method yields
a better error estimate in the ASIP than the usual one,
see Field~{\em et al.}~\cite{FMT03}.
However, it should be emphasised that the martingale approximation
of Gordin~\cite{Gordin69} leads directly only to
 a {\em reverse} martingale increment sequence
and so the ASIP is obtained in backwards time in the first instance.
This is not an issue for distributional results such as the CLT,
but the ASIP in~\cite{ConzeBorgne01,FMT03,MT02} uses explicitly the fact
that the class of systems being studied is closed under time-reversal.
To obtain {\em forward} martingale approximations, it is necessary
to use more sophisticated versions of Gordin's approach~\cite{PhilippStout75}.

Recently, there has been an explosion of interest in nonuniformly
expanding maps and nonuniformly hyperbolic diffeomorphisms
(possibly with singularities).  We refer to the articles
of Young~\cite{Young98,Young99} as well as Aaronson~\cite{Aaronson},
Baladi~\cite{Baladi,Baladi01},
Gou\"ezel~\cite{GouezelPhD}, Viana~\cite{Viana} and references therein.
In particular, decay of correlations and the CLT are studied extensively
in these references.
However, such classes of dynamical systems are intrinsically time-orientation
specific, and largely for this reason the ASIP has not previously been proved.
Similarly, the LIL was previously unproved for such systems.

In this paper, we establish the ASIP, and hence the (functional) LIL,
 for nonuniformly expanding/hyperbolic systems.   Both discrete time systems
and flows are covered by our results.

\begin{rmk} \label{rmk-PS}
We note that~\cite{PollicottSharp02} attempted to apply 
the approach in~\cite{FMT03} to nonuniformly expanding systems.
However, it appears that the time-orientation issue discussed above
was overlooked in~\cite{PollicottSharp02}, and that this is a gap.
Hence it seems necessary to find an alternative approach to the one in~\cite{FMT03},
and that is what is done in the current paper.
\end{rmk}

Precise formulations are given in the body of the paper, but here 
is an outline of our main result, and the strategy behind its proof, for a 
nonuniformly expanding map $T:M\to M$ where $(M,d)$ is a metric space.
By standard methods, $T$ can 
be modelled by a discrete-time suspension over a {\em Gibbs-Markov}
map~\cite{Aaronson} $f:Y\to Y$ with return time function $R:Y\to\Z^+$.   
(Roughly speaking, a Gibbs-Markov map is like a uniformly expanding map
with possibly countably many inverse branches.)
There exists a unique ergodic $T$-invariant probability measure 
equivalent to Lebesgue, and the following result is formulated with
this measure in mind.

\begin{thm} \label{thm-intro}   Let $T:M\to M$ be a nonuniformly
expanding map.  Assume moreover that $R\in L^{2+\delta}(Y)$.
Let $\phi:M\to\R$ be a mean zero H\"older observation.
Then $\phi$ satisfies the ASIP.
That is, there exists $\epsilon>0$, a sequence of random variables $\{S_n\}$
and a Brownian motion $W$ with variance $\sigma^2\ge0$ 
such that $\{\sum_{j=0}^{N-1}\phi\circ T^j\}=_d\{S_N\}$, and
\[
S_N=W(N)+O(N^{\frac12-\epsilon}) \quad\text{as $N\to\infty$},
\]
almost everywhere.
\end{thm}

Using a method due to Hofbauer and Keller~\cite{HofbauerKeller82}
which exploits a result of Philipp and Stout~\cite[Theorem~7.1]{PhilippStout75},
we obtain the ASIP (in the correct time direction but without the improved
error term) for $f:Y\to Y$ and a class of ``weighted Lipschitz'' 
observations.   Theorem~\ref{thm-intro} then
follows directly by Melbourne and T\"or\"ok~\cite{MT04}.
(We note that the method in~\cite{MT04} has independently been used by 
Gou\"ezel~\cite{Gouezelsub_a}
to obtain a simplified derivation of the CLT and stable laws.)

A precise version of Theorem~\ref{thm-intro} is stated and proved in
Section~\ref{sec-expand}(e).
The ASIP for nonuniformly hyperbolic maps extends easily to a class
of nonuniformly expanding semiflows, see Section~\ref{sec-expand}(e).

Our results for nonuniformly hyperbolic diffeomorphisms
and nonuniformly hyperbolic flows are completely analogous, but the set-up
is more technical and we postpone further details until Section~\ref{sec-hyp}.

\vspace{-2ex}
\paragraph{Planar periodic Lorentz gas}
The planar periodic Lorentz gas is a class of examples
introduced by Sina{\u\i}~\cite{Sinai70}.
See~\cite{ChernovYoung00} for a survey
of results about Lorentz gases.   The Lorentz flow is a billiard flow on
$\T^2-\Omega$
where $\Omega$ is a disjoint union of convex regions with $C^3$ boundaries.
(The phase-space of the flow is three-dimensional; planar position and
direction.)   The flow has a natural global cross-section
$M=\partial\Omega\times[-\pi/2,\pi/2]$ corresponding to collisions
and the Poincar\'e map $T:M\to M$ is called the billiard map.
Bunimovich, Sina{\u\i} and Chernov~\cite{BunimSinaiChernov91} proved
the central limit theorem and weak invariance principle for such maps.

Denote the return time function by $h:M\to\R^+$.  The Lorentz flow
satisfies the {\em finite horizon} condition if $h$ is uniformly bounded.
The central limit theorem and weak invariance principle was proved 
by~\cite{BunimSinaiChernov91} for Lorentz flows satisfying
the finite horizon condition.

\begin{thm} \label{thm-Lorentz}
Suppose that $T_t$ is a planar periodic Lorentz gas.
\begin{itemize}
\item[(i)] The billiard map satisfies the ASIP for H\"older observations.
\item[(ii)] If the finite horizon condition holds, then the Lorentz flow
satisfies the ASIP for H\"older observations.
\end{itemize}
\end{thm}

In Section~\ref{sec-expand}, we prove the 
ASIP for nonuniformly expanding maps and semiflows.
In Section~\ref{sec-hyp}, we prove the analogous results for systems that
are nonuniformly hyperbolic in the sense of Young~\cite{Young98}.
In Section~\ref{sec-app}, we list numerous examples in the literature for 
which our results apply.  In particular, we prove Theorem~\ref{thm-Lorentz}.
The results of~\cite{PhilippStout75,MT04} required in this paper
are reproduced as appendices.

\section{Nonuniformly expanding systems}
\label{sec-expand}

In this section, we prove the ASIP for nonuniformly expanding systems.
The first step is to prove the ASIP for Gibbs-Markov maps.
Such maps are reviewed in Subsection~(a) and a class of ``weighted
Lipschitz'' observations is introduced in Subsection~(b).
The ASIP for Gibbs-Markov maps is proved in Subsection~(c)
using an approach of Hofbauer and Keller~\cite{HofbauerKeller82}.
In Subsection~(d), we obtain the ASIP for
Young towers~\cite{Young99} as an application of~\cite{MT04}.   In
Subsection~(e), we prove the ASIP for nonuniformly expanding maps and 
semiflows.

\subsection{Gibbs-Markov maps}

Let $(\Lambda,m)$ be a Lebesgue space with a countable
measurable partition $\alpha$.
Without loss, we suppose that all partition elements $a\in\alpha$ have 
$m(a)>0$.
Recall that a measure-preserving transformation $f:\Lambda\to \Lambda$ is a 
{\em Markov map} if $f(a)$ is a union of elements  of $\alpha$
and $f|_a$ is injective for all $a\in\alpha$.
Define $\alpha'$ to be the coarsest partition of $\Lambda$ such that 
$fa$ is a union of atoms in $\alpha'$ for all $a\in\alpha$.
(So $\alpha'$ is a coarser partition than $\alpha$.)
If $a_0,\ldots,a_{n-1}\in\alpha$, we define the £$n$-cylinder
$[a_0,\ldots,a_{n-1}]=\cap_{i=0}^{n-1}f^{-i}a_i$.
It is assumed that $f$ and $\alpha$ separate points in $\Lambda$
(if $x,y\in \Lambda$ and $x\neq y$, then for $n$ large enough there
exist distinct $n$-cylinders that contain $x$ and $y$).

Let $0<\beta<1$.
We define a metric $d_\beta$ on $\Lambda$ by $d_\beta(x,y)=\beta^{s(x,y)}$ 
where $s(x,y)$ is the greatest integer $n\ge0$ such that
$x,y$ lie in the same $n$-cylinder.
Define $g=Jf^{-1}=\frac{dm}{d(m\circ f)}$ and set
$g_k=g\,g\circ f\,\cdots\,g\circ f^{k-1}$.

The map $f:\Lambda\to \Lambda$ is a {\em Gibbs-Markov map} if it satisfies 
the additional 
properties:
\begin{itemize}
\item[(i)] {\em Big images property:}   There exists $c>0$ such that
$m(fa)\ge c$ for all $a\in\alpha$.
\item[(ii)] {\em Distortion:}  $\log g|_a$ is Lipschitz with respect to
$d_\beta$ for all $a\in\alpha'$.
\end{itemize}
It follows from assumptions~(i) and~(ii) that there exists a 
constant $D\ge1$ such that for all $x,y$ lying in a common $k$-cylinder
$[a_0,\ldots,a_{k-1}]$,
\begin{align} \label{eq-distort}
\Bigl|\frac{g_k(x)}{g_k(y)}-1\Bigr|\le Dd_\beta(f^kx,f^ky)
\quad\text{and}\quad  
D^{-1}\le \frac{m[a_0,\ldots,a_{k-1}]}{g_k(x)} \le D.
\end{align}

\subsection{Weighted Lipschitz observations}

Let $p\ge1$.  We fix a sequence of weights $R(a)>0$ satisfying
$|R|_p=(\sum_{a\in\alpha} m(a) R(a)^p)^{1/p} < \infty$.

Given $v:\Lambda\to\R$ continuous, we set $v_a=v|_a$
and define $|v|_\beta$ to be the Lipschitz constant of $v$ with respect
to the metric $d_\beta$.  Let
\[
\|v\|_\infty = \sup_{a\in\alpha} |v_a|_\infty/R(a), \qquad 
\|v\|_\beta=\sup_{a\in\alpha} |v_a|_\beta/R(a).
\]
Let ${\cal B}$ consist of the space of weighted Lipschitz functions
with \mbox{$\|v\|=\|v\|_\infty +\|v\|_\beta<\infty$}.
Note in particular that $R\in {\cal B}$ and $\|R\|=1$.
We have the embeddings
\[
{\rm Lip}\subset {\cal B}\subset L^p \subset L^1,
\]
where ${\rm Lip}$ is the space of (globally) Lipschitz functions.

The transfer (Perron-Frobenius) operator $P:L^1\to L^1$ 
maps $v\in L^1$ to $Pv$ where 
$\int_\Lambda Pv\,w\,dm=\int_\Lambda v\,w\circ f\,dm$
for all $w\in L^\infty$, and is given by $(Pv)(x)=\sum_{fy=x}g(y)v(y)$. 
Note that $|P|_1=1$.

\begin{prop} \label{prop-basic}
Let $a\in\alpha$ be an $n$-cylinder and suppose that $v:a\to\R$
is Lipschitz.   Then
$|v|_\infty \le \frac{1}{m(a)}\int_a |v|\,dm+\beta^n|v|_\beta$.
\end{prop}

\begin{proof}
For $x\in a$,
\begin{align*}
|v(x)| &\le  \SMALL\frac{1}{m(a)}\int_a |v|\,dm + |v(x)- \frac{1}{m(a)}\int_a v\,dm|
\le \frac{1}{m(a)}\int_a |v|\,dm+|v|_\beta\diam(a).
\end{align*}
The result follows since $\diam(a)=\beta^n$.
\end{proof}

\begin{lemma} \label{lem-basic}
The transfer operator $P$ restricts to an operator $P:{\cal B}\to{\cal B}$
and there exists a constant $C\ge1$ such that
\[
\|P^nv\|\le C(|v|_1 + \beta^n \|v\|_\beta),
\]
for all $v\in\mathcal{B}$ and $n\ge 1$.
Moreover, $P({\cal B})\subset {\rm Lip}$.
\end{lemma}

\begin{proof}
We prove the estimate on $\|P^nv\|$.
The remaining statements of the lemma are evident from the proof.

Note that $(P^n v)(x)=\sum_{f^ny=x}g_n(y)v(y)$.
Since the $n$-cylinders $[a_0,\ldots,a_{n-1}]$ form a partition
and each $n$-cylinder contains precisely one preimage $y_a$, we have
\begin{align*}
|(P^nv)(x)| & \le \sum_{a=[a_0,\ldots,a_{n-1}]}g_n(y_a)|v(y_a)| 
\le D\sum_{a=[a_0,\ldots,a_{n-1}]}m(a)|v_a|_\infty \\
& \le D\sum_{a=[a_0,\ldots,a_{n-1}]}\Bigl[{\SMALL \int}_a|v|\,dm+m(a)\beta^n|v_a|_\beta\Bigr] \\
& \le D\sum_{a=[a_0,\ldots,a_{n-1}]}\Bigl[{\SMALL \int}_a|v|\,dm+\beta^nm(a)R(a_0)\|v\|_\beta\Bigr] 
\end{align*}
where we have used Proposition~\ref{prop-basic} and estimate~\eqref{eq-distort}.
Hence $|P^nv|_\infty \le D\Bigl[|v|_1+\beta^n|R|_1\|v\|_\beta\Bigr]$.
Similarly,
\begin{align*}
|(P^nv)(x)- (P^nv)(x')| & 
\le \sum_{a=[a_0,\ldots,a_{n-1}]}|g_n(y_a)-g_n(y_a')|
|v(y_a)| \\ &+
\sum_{a=[a_0,\ldots,a_{n-1}]}|g_n(y_a')|
|v(y_a)- v(y_a')|.
\end{align*}
Each term in the first summation can be estimated by
\begin{align*}
D|g_n(y_a')|d_\beta(f^ny_a,f^ny_a')|v_a|_\infty
& \le D^2 m(a)d_\beta(x,x')
\bigl[{\SMALL\frac{1}{m(a)}\int_a |v_a|\,dm+\beta^n|v_a|_\beta}\bigr] \\
&\le D^2\Bigl[{\SMALL \int}_a|v|\,dm+\beta^n m(a)R(a_0)\|v\|_\beta\Bigr]d_\beta(x,x),
\end{align*}
so the first summation is bounded by
$D^2\bigl[|v|_1+\beta^n|R|_1\|v\|_\beta\bigr] d_\beta(x,x')$.
Each term in the second summation can be estimated by
\begin{align*}
Dm(a)|v_a|_\beta d_\beta(y_a,y_a')
\le Dm(a)R(a_0)\|v\|_\beta \beta^n d_\beta(x,x'),
\end{align*}
so the second summation is bounded by
$\beta^n D|R|_1\|v\|_\beta d_\beta(x,x')$.
The result follows.
\end{proof}

We have the following standard consequences of Lemma~\ref{lem-basic}.

\begin{cor} \label{cor-qc}
Let $p,q\ge1$ with $\frac1p+\frac1q=1$.
Assume that $f:\Lambda\to\Lambda$ is mixing and that $R\in L^p$.   Then
there exist constants $C\ge1$ and $\tau\in(0,1)$ such that
\begin{itemize}
\item[(a)]
$\|P^nv-\int_\Lambda v\,dm\|\le C\tau^n\|v\|$ for all 
$v\in\mathcal{B}$ and $n\ge1$.
\item[(b)]  $|\int_\Lambda v\,(w\circ f^n)\,dm - \int_\Lambda v\,dm
\int_\Lambda w\,dm| \le C\tau^n\|v\| |w|_q$
for all $v\in\mathcal{B}$, $w\in L^q$, $n\ge1$.
\item[(c)] If $R\in L^2$, then for any $v\in\mathcal{B}$ with $\int_\Lambda v\,dm=0$, 
the series
\[
\sigma^2=\int_\Lambda v^2\,dm+2\sum_{k=1}^\infty \int_\Lambda v\,(v\circ f^k)\,dm,
\]
is absolutely convergent, and
$\int_\Lambda v_N^2\,dm=\sigma^2N+O(1)$ as $N\to\infty$,
where $v_N=\sum_{j=0}^{N-1}v\circ f^j$.
Moreover,
$\sigma=0$ if and only if there exists a Lipschitz function $w:\Lambda\to\R$
such that $v=w\circ f-w$.
\end{itemize}
\end{cor}

\begin{proof}   Most of this result is completely standard, but we include the
details for completeness.  By an Arzela-Ascoli argument, the unit ball
in $\mathcal B$ is compact in $L^1$.   This combined with
Lemma~\ref{lem-basic} implies, by Hennion~\cite{Hennion93},
that the essential spectral radius of $P:\mathcal{B}\to\mathcal{B}$
is bounded above by $\beta<1$.   There is a simple eigenvalue at $1$
with eigenspace consisting of constant functions, but the mixing assumption
guarantees that there are no further eigenvalues on the unit circle.
Now choose $\tau\in(\beta,1)$ such that all eigenvalues of $P$
other than $1$ lie strictly inside the disk of radius $\tau$.
Part~(a) follows for such a choice of $\tau$.

To prove part~(b), compute that 
\begin{align*}
|{\SMALL \int}_\Lambda v\,(w\circ f^n)\,dm - {\SMALL \int}_\Lambda v\,dm
{\SMALL \int}_\Lambda w\,dm| 
& = |{\SMALL \int}_\Lambda (P^nv-{\SMALL \int}_\Lambda v)\,w\,dm|
\le |P^nv-{\SMALL \int}_\Lambda v|_p |w|_q \\
& \le \|P^nv-{\SMALL \int}_\Lambda v\| |w|_q \le C\tau^n\|v\| |w|_q.
\end{align*}

It follows from (b) that $|{\SMALL \int}_\Lambda v\,(v\circ f^k)\,dm|\le C\tau^n\|v\| |v|_2$
and so the series for $\sigma^2$
converges absolutely.   Moreover
\begin{align*}
{\SMALL \int}_\Lambda v_N^2\,dm & = N{\SMALL \int}_\Lambda v^2\,dm+2\sum_{0\le i<j\le N-1}
{\SMALL \int}_\Lambda v\,(v\circ f^{j-i})\,dm  \\ &=
N{\SMALL \int}_\Lambda v^2\,dm+2\sum_{k=1}^N (N-k){\SMALL \int}_\Lambda v\,(v\circ f^k)\,dm \\
& = N\sigma^2 - 2\sum_{k=1}^N k{\SMALL \int}_\Lambda v\,(v\circ f^k)\,dm
-2\sum_{k=N+1}^\infty N{\SMALL \int}_\Lambda v\,(v\circ f^k)\,dm \\
&= N\sigma^2+O(1),
\end{align*}
proving (c).

The criterion for $\sigma=0$ follows as in~\cite{FMT03,MNapp}.
If $v=w\circ f-w$, then $v_N=w\circ f^N-w$ so it is clear that $\sigma=0$.
To prove the converse,
define $w=\sum_{j=1}^\infty P^jv$.   This series converges in $\mathcal{B}$
by~(b) and is Lipschitz by Lemma~\ref{lem-basic}.
Write $v=\hat v+w\circ f-w$.   Then it is easily
seen that $\hat v$ has the same variance as $v$ and that $P\hat v=0$.
Hence $\sigma^2={\SMALL \int}_\Lambda \hat v^2\,dm$, so if $\sigma=0$, then 
$\hat v=0$ and $v=w\circ f-w$.
\end{proof}

\subsection{ASIP for Gibbs-Markov maps}

Let $\alpha_0^{k-1}$ denote the partition into length $k$ cylinders
$a=[a_0,\ldots,a_{k-1}]$.

\begin{lemma} \label{lem-PS}
Assume that $f:\Lambda\to\Lambda$ is mixing and that
$R\in L^{2+\delta}$ for some $\delta>0$.
Let $v\in\mathcal{B}$ with $\int_\Lambda v\,dm=0$.
Then
\begin{itemize}
\item[(a)]  $\sum_{a\in\alpha_0^{k-1}} 
\int_a|v-\frac{1}{m(a)}\int_a v\,dm|^{2+\delta}\,dm \le \bigl(\|v\|_\beta\,|R|_{2+\delta}\,\beta^k \bigr)^{2+\delta}$.
\item[(b)]  
$\bigl|m(a\cap f^{-(N+k)}(b))-m(a)m(b)\bigr|\le C\tau^Nm(a)m(b)^{1/2}$
for all $a\in\alpha_0^{k-1}$ and all measurable sets $b$.
\end{itemize}
\end{lemma}

\begin{proof}
Note that
$|v-\frac{1}{m(a)}\int_a v\,dm|\le |v_a|_\beta\,\diam(a)\le \|v\|_\beta
R(a_0)\beta^k$.
Part~(a) follows immediately.

We argue as in Aaronson \& Denker~\cite{AaronsonDenker01} to establish (b).
Let $v_{a,k}=P^k\chi_a$.
By definition, $v_{a,k}=\sum_{f^ky=x}g_k(y)\chi_a(y)=g_k(y_a)$
where $y_a$ is the unique point in $a$ such that $f^ky_a=x$.
Hence by~\eqref{eq-distort},
\[
|v_{a,k}(x)-v_{a,k}(x')|\le D|g_k(y_a)|d_\beta(x,x')\le D^2m(a)d_\beta(x,x').
\]
It follows that $\|v_{a,k}\|\le Em(a)$ where $E=D^2+D$.

Using this estimate and Corollary~\ref{cor-qc}(b), we compute that
\begin{align*}
& \bigl|m(a\cap f^{-(N+k)}b)-m(a)m(b)\bigr| =
\bigl|{\SMALL \int} P^k\chi_a\,\chi_b\circ f^N-{\SMALL \int} P^k\chi_a{\SMALL \int} \chi_b\bigr|  \\
& \qquad = 
\bigl|{\SMALL \int} v_{a,k}\,\chi_b\circ f^N-{\SMALL \int} v_{a,k}
{\SMALL \int} \chi_b\bigr| \le C\tau^N\|v_{a,k}\|\,|\chi_b|_2  
\le CE\tau^N m(a)m(b)^{1/2},
\end{align*}
as required.
\end{proof}

\begin{cor}   \label{cor-ASIPGM}
Let $f:\Lambda\to\Lambda$ be an ergodic Gibbs-Markov map.
Define the Banach space $\mathcal{B}$ corresponding to 
weights $R\in L^{2+\delta}$ for some $\delta>0$.
Suppose that $v\in\mathcal{B}$ and $\int_\Lambda v\,dm=0$.
Define $\sigma^2$ as in Corollary~\ref{cor-qc} and assume that $\sigma^2>0$.
Then $v_N=\sum_{j=0}^{N-1}v\circ f^j$ satisfies the ASIP.
\end{cor}

\begin{proof}
We verify the hypotheses of Philipp \& Stout~\cite[Theorem~7.1]{PhilippStout75}.
For convenience, we have translated this theorem into dynamical systems
terminology in the appendix, see Theorem~\ref{thm-PS}.
Condition~(i) of Theorem~\ref{thm-PS} is automatic since 
$\mathcal{B}\subset L^{2+\delta}$
and condition~(ii) follows from Corollary~\ref{cor-qc}(c).
Conditions~(iii) and~(iv) follow from parts~(a) and~(b)
of Lemma~\ref{lem-PS}.
\end{proof}

\subsection{ASIP for tower maps}

Suppose that $(\Lambda,m)$ is a probability space and that
$f:\Lambda\to\Lambda$ is a measure-preserving transformation.
Let $R:\Lambda\to\Z^+$ be a measurable function (called a
{\em return time function} with $R\in L^1(\Lambda)$.
Define the suspension 
\begin{align*}
\Delta=\{(x,\ell)\in \Lambda\times\N: 0\le \ell\le R(x)\}/\sim,
\end{align*}
where $(x,R(x))\sim(f(x),0)$.
Define $F:\Delta\to\Delta$ by $F(x,\ell)=(x,\ell+1)$ computed
subject to identifications.
Note in particular that $F(x,0)=(f(x),0)$.
An $F$-invariant probability measure on $\Delta$ is given by
$m^R=m\times l/\bar R$ where $\bar R=\int_{\Lambda} R\,dm$
and $l$ is counting measure on $\N$.  

Let $\{\Delta_{j,0}\}$ be a countable measurable partition of $\Lambda$ 
such that $f$ and $\{\Delta_{j,0}\}$ separate points in 
$\Lambda$, and for each $j$, $R_j=R|_{\Delta_{j,0}}$ is constant and
 $f:\Delta_{j,0}\to\Lambda$ is a measurable isomorphism.
For each $j$ and $0\le \ell<R_j$, let 
$\Delta_{j,\ell}=\Delta_{j,0}\times\{\ell\}$.
This defines a partition $\{\Delta_{j,\ell}\}$ of $\Delta$.

A {\em separation time} function $s:\Delta\times\Delta\to\N$ is defined
as follows:   If $x,y$ lie in distinct partition elements, then $s(x,y)=0$.
If $x,y\in\Delta_{j,0}$ for some $j$, then $s(x,y)$ is the greatest
integer $n\ge0$ such that $f^kx$ and $f^ky$ lie in the same partition
element of $\Lambda$ for $k=0,\ldots,n$.   
If $x,y\in\Delta_{j,\ell}$, then write $x=F^\ell x_0$,
$y=F^\ell y_0$ where $x_0,y_0\in\Delta_{j,0}$ and define $s(x,y)=s(x_0,y_0)$.
For $\theta\in(0,1)$, we define a metric $d_\theta$ on $\Delta$
by setting $d_\theta(x,y)=\theta^{s(x,y)}$.

\begin{defn}  The suspension $F:\Delta\to\Delta$ is called a {\em Young tower}
if $f:\Lambda\to\Lambda$ is a Gibbs-Markov map with respect
to the partition $\alpha=\{\Delta_{j,0}\}$.
\end{defn}

\begin{rmk}  The big images condition for $f$ to be a Gibbs-Markov map
is automatically satisfied in the strong sense that $f(a)=\Lambda$
for each each $a\in\alpha$.   Hence, $F:\Delta\to\Delta$ is a Young tower
provided the distortion condition holds: there exist constants $\theta\in(0,1)$
and $C\ge1$ such that for each $j$ the Jacobian $g_j=Jf|_{\Delta_{j,0}}:\Delta_{j,0}\to\Lambda$ satisfies $|\log g_j(x)-\log g_j(y)|\le Cd_\theta(x,y)$
for all $x,y\in\Delta_{j,0}$.
\end{rmk}

\begin{thm} \label{thm-ASIPtower}   Let $F:\Delta\to\Delta$ be a Young tower
defined as a suspension over $f:\Lambda\to\Lambda$ with return time function $R$.
Assume that $R\in L^{2+\delta}(\Lambda)$.
Let $\phi:\Delta\to\R$ be a mean zero observation and assume
that $\phi$ is Lipschitz with respect to $d_\theta$.
Then $\phi_N=\sum_{j=0}^{N-1}\phi\circ F^j$ satisfies the ASIP.
\end{thm}

\begin{proof}
Define a mean zero observation $\Phi:\Lambda\to\R$ by setting
$\Phi(x)=\sum_{j=0}^{R(x)-1}\phi(x,j)$.
Since $\phi$ is Lipschitz, it is immediate that $\Phi$ lies in
the space $\mathcal{B}$ of weighted Lipschitz observations.
Since $R\in L^{2+\delta}(\Lambda)$, it follows from Corollary~\ref{cor-ASIPGM}
that $\Phi_N=\sum_{j=0}^{N-1}\Phi\circ f^j$ satisfies the ASIP on $\Lambda$.

Note that $R-\bar R$ also satisfies the hypotheses of 
Corollary~\ref{cor-ASIPGM}, and so the ASIP, and hence the LIL, applies.   
Therefore, it is certainly the case
that $\sum_{j=0}^{N-1}R\circ f^j=N\bar R + o(N^{1-\delta})$ almost everywhere.
The result follows from~\cite[Theorem~4.2]{MT04}, see Corollary~\ref{cor-MT}.
\end{proof}

\subsection{ASIP for nonuniformly expanding systems}

Let $(M,d)$ be a locally compact separable bounded
metric space with Borel probability
measure $\eta$ and let $T:M\to M$ be a nonsingular transformation for
which $\eta$ is ergodic.
Let $Y\subset M$ be a measurable subset with $\eta(Y)>0$.
We suppose that there is an at
most countable measurable partition $\{Y_j\}$ with $\eta(Y_j)>0$,
and that there exist integers $R_j\ge1$, and constants $\lambda>1$;
$C,D>0$ and $\gamma\in(0,1)$ such that for all $j$,
\begin{itemize}
\item[(1)]  $T^{R_j}:Y_j\to Y$ is a (measure-theoretic) bijection.
\item[(2)]  $d(T^{R_j}x,T^{R_j}y)\ge \lambda d(x,y)$ for all $x,y\in Y_j$.
\item[(3)] $d(T^kx,T^ky)\le Cd(T^{R_j}x,T^{R_j}y)$ for all $x,y\in Y_j$,
$k<R_j$.
\item[(4)] $g_j=\frac{d(\eta|_{Y_j}\circ (T^{R_j})^{-1})}{d\eta|_Y}$ satisfies
$|\log g_j(x)-\log g_j(y)|\le Dd(x,y)^\gamma$ for almost all $x,y\in Y$.
\item[(5)] $\sum_j R_j\eta(Y_j)<\infty$.
\end{itemize}
We say that a dynamical system $T$ satisfying (1)--(5) is
{\em nonuniformly expanding}.

Define the {\em return time function} $R:Y\to \Z^+$ by $R|_{Y_j}\equiv R_j$.
Condition~(5) says that $\int_Y R\,d\eta<\infty$.
The map $f:Y\to Y$ given by $f(y)=T^{R(y)}(y)$
is the corresponding {\em induced map}.
It can be shown (see Young~\cite[Theorem~1]{Young99}) that
there is a unique invariant
probability measure $m$ on $M$ that is equivalent to $\eta$.

We can now state and prove a precise version of Theorem~\ref{thm-intro}.

\begin{thm} \label{thm-expand}   Let $T:M\to M$ be a nonuniformly
expanding map satisfying (1)--(5) above.
Assume moreover that the return time function $R$ lies in $L^{2+\delta}(Y)$.
Let $\phi:M\to\R$ be a mean zero H\"older observation.
Then $\phi$ satisfies the ASIP.
\end{thm}

\begin{proof}
Let $\Delta=\{(y,\ell):y\in Y,\,\ell=0,\dots,R(y)-1\}$,
so $\Delta$ is the disjoint union of $R_j$ copies of each $Y_j$.
Define a measure $\mu$ on $\Delta$ by setting
$\mu|_{Y_j\times\{\ell\}}=m|_{Y_j}/\bar R$.
Define $F:\Delta\to \Delta$ by setting 
$F(y,\ell)=(y,\ell+1)$ for $0\le \ell<R(y)-1$
and $F(y,R(y)-1)=(fy,0)$. Define  the separation time
$s:\Delta \times \Delta \to \N$ as in the previous section.

By shrinking $\gamma$ if necessary, we may suppose that $\phi$ is 
$\gamma$-H\"older for the same $\gamma$ that appears in condition~(4).
Define the metric $d_\theta$ on $\Delta$ with $\theta=1/\lambda^\gamma$.
It follows from condition~(2) 
that $d(x,y)\le \diam(Y)/\lambda^{s(x,y)}$ for all $(x,y)\in\Delta$.
Hence $f$ and $\{Y_j\}$ separate points in $Y$
and the required distortion condition on $g_j$ is immediate, so
$\Delta$ is a Young tower with $\Lambda=Y$ and $\Delta_{j,0}=Y_j$.

If $x,y$ lie in the same partition element of $\Delta_\ell$, then
write $x=F^\ell x_0$, $y=F^\ell y_0$ so
$d(fx_0,fy_0)\le \diam(Y)/\lambda^{s(x,y)}$.
By condition~(3),
\begin{align*}
d(x,y)\le C d(fx_0,fy_0) \le C \diam(Y)/\lambda^{s(x,y)}
= C\diam(Y)[d_\theta(x,y)]^{1/\gamma}.
\end{align*}
Hence, there is a constant $C'\ge1$ such that
$d(x,y)\le C'[d_\theta(x,y)]^{1/\gamma}$ for all $x,y\in\Delta$.

Define the projection $\pi:\Delta\to M$ by $\pi(y,\ell)=T^\ell y$.
Then $\pi$ is a measure-preserving isomorphism and it follows as above that
$d(\pi(x),\pi(y))^\gamma\le C''d_\theta(x,y)$, for all $x,y\in\Delta$.
In particular, since $\phi:(M,d)\to\R$ is $\gamma$-H\"older,
it follows that $\phi\circ\pi:(\Delta,d_\theta)\to\R$ is Lipschitz.
By Theorem~\ref{thm-ASIPtower}, the ASIP holds for $\phi\circ\pi$
on $\Delta$.    Since $\pi$ is a measure-preserving map semiconjugacy, the ASIP
holds for $\phi$ on $M$.
\end{proof}

\begin{rmk} \label{rmk-CLT}
As already pointed out in~\cite{Gouezelsub_a},
the CLT for nonuniformly expanding maps holds under slightly weaker 
hypotheses using~\cite[Theorem~1.1]{MT04}.   Instead of requiring that 
$R\in L^{2+\delta}$, it suffices that $R\in L^2$.
\end{rmk}

\begin{rmk} \label{rmk-deg}
The ASIP is said to be {\em degenerate} if $\sigma^2=0$.
It follows from previous work in connection with the CLT~\cite{Young98,Young99}
that the ASIPs obtained in this paper are degenerate if and only if 
$\phi=\psi\circ T-\psi$ where
$\psi\in L^2(M)$.  Moreover, by a Liv\v{s}ic regularity result of
Bruin {\em et al.}~\cite{BHN}, such an $L^2$ function $\psi$ has a version
that is H\"older on $\cup_{j=0}^\ell T^jY$ for each fixed~$\ell$.
(It is easy to construct examples where the ASIP
is degenerate but $\psi$ does not have a version that is continuous on the
whole of $M$.)  
In particular, if $T$ has a periodic point $x\in Y$ of period $k$
and $\sum_{j=0}^{k-1}\phi(T^jy)\neq0$, then the ASIP 
is nondegenerate.
\end{rmk}

\paragraph{Nonuniformly expanding semiflows}
We continue to assume that $T:M\to M$ is a nonuniformly expanding map
satisfying conditions (1)--(5).
Suppose that $h:M\to\R^+$ lies in $L^1(M)$.
Regarding $h$ as a {\em roof function}, we form the {\em suspension}
$M^h=\{(x,u)\in M\times[0,\infty):0\le u\le h(x)\}/\sim$
where $(x,h(x))\sim(Tx,0)$.   The suspension semiflow $T_t:M^h\to M^h$
is given by $T_t(x,u)=x(u+t)$ computed modulo identifications.
We call $T_t:M^h\to M^h$ a {\em nonuniformly expanding semiflow}.
We say that an observation $\psi:M^h\to \R$ is H\"older if $\psi$ is bounded and
$\sup_{(x,u)\neq(y,u)}|\psi(x,u)-\psi(y,u)|/d(x,y)<\infty$.

\begin{cor} \label{cor-semiflow} Let $T_t:M^h\to M^h$ be a nonuniformly
expanding semiflow.   Assume moreover that 
the return time function $R$ lies in $L^{2+\delta}(Y)$ and that
the roof function $h:M\to\R^+$ is H\"older.
Let $\psi:M^h\to\R$ be a mean zero H\"older observation.
Then $\psi$ satisfies the ASIP.
That is, there exists $\epsilon>0$, a family of random variables $\{S_t\}$
and a Brownian motion $W$ with variance $\sigma^2\ge0$ 
such that $\{\int_0^t \psi\circ T_s\,ds\}=_d\{S_t\}$, and
$S_t=W(t)+O(t^{\frac12-\epsilon})$ as $t\to\infty$, almost everywhere.
\end{cor}

\begin{proof}
According to~\cite[Theorem~4.2]{MT04} (Theorem~\ref{thm-MT}), it suffices that (i) 
$h\in L^{2+\delta}(Y)$, (ii) $\phi(x)=\int_0^{h(x)}\psi(x,u)du$
satisfies the ASIP on $Y$, and (iii) $h$ satisfies the ASIP on $Y$.
Hence, the result is immediate from Theorem~\ref{thm-expand}.
\end{proof}

\begin{rmk} \label{rmk-semiflow}
We have not striven for greatest generality
in the statements of Theorem~\ref{thm-expand} and
Corollary~\ref{cor-semiflow}. 
However, it is clear from the proof that in Theorem~\ref{thm-expand}
we can relax the assumption that $\phi$ is H\"older.  It is sufficient
that $\phi$ is such that $\Phi(x)=\sum_{\ell=0}^{R(x)-1}\phi(T^\ell x)$
lies in the space of weighted Lipschitz observations in Subsection~(b)
for an appropriate choice of weight function.   Taking the weight
function to be the return time function, it suffices that
$\phi$ is H\"older on $T^\ell Y_j$ for all $j\ge1$, $0\le \ell<R(j)-1$,
with $L^\infty$ norm and H\"older constant independent of $j,\ell$.

Similarly, the hypotheses that $\psi$ and $h$ are
H\"older can be weakened in  Corollary~\ref{cor-semiflow}.
For example, provided $\psi$ is H\"older, it suffices that $h$ 
is H\"older on $T^\ell Y_j$ for all $j\ge1$, $0\le \ell<R(j)-1$,
with $L^\infty$ norm and H\"older constant independent of $j,\ell$.
\end{rmk}

\section{Nonuniformly hyperbolic systems}
\label{sec-hyp}

In this section, we show how to prove the ASIP for Lipschitz observations
of a dynamical system that is {\em nonuniformly hyperbolic
in the sense of Young~\cite{Young98}}.
Instead of using the original set up, we make four assumptions
(A1)--(A4) that are distilled from those in~\cite{Young98}.
In doing so, we bypass the differential structure, and certain conclusions
in~\cite{Young98} become assumptions here, particularly (A4) below.

Let $T:M\to M$ be a diffeomorphism (possibly with singularities) defined on a 
Riemannian manifold $(M,d)$.   We assume from the start that $T$ preserves
a ``nice'' probability measure $m$ (one of
the conclusions in Young~\cite{Young98} is that $m$ is a SRB measure).
Assumption (A4) contains the properties of $m$ that we require for the ASIP.

We fix a subset $\Lambda\subset M$ and a family of subsets of $M$ that
we call ``stable 
disks'' $\{W^s\}$ that are disjoint and cover $\Lambda$.
If $x$ lies in a stable disk, we label the disk $W^s(x)$.

\begin{itemize}
\item[(A1)]  There is a partition $\{\Lambda_j\}$ of $\Lambda$ and
integers $R_j\ge1$ such that 
for all $x\in \Lambda_j$ we have $T^{R_j}(W^s(x))\subset W^s(T^{R_j}x)$.
\end{itemize}

Define the return time function $R:\Lambda\to\Z^+$ by $R|_{\Lambda_j}=R_j$
and the induced map
$f:\Lambda \to \Lambda$ by $f(x)=T^{R(x)} (x)$.
Form the discrete suspension map $F:\Delta\to\Delta$ where 
$F(x,\ell)=(x,\ell+1)$ for $\ell<R(x)-1$ and $F(x,R(x)-1)=(fx,0)$.
We define a separation time  $s:\Lambda\times \Lambda \to \N$ by defining 
$s(x,x')$ to be the greatest integer $n\ge0$
such that $f^k x, f^k x'$ lie in the same partition element of $\Lambda$
for $k=0,\ldots,n$.
(If $x,x'$ do not lie in the same partition element, then we take
$s(x,x')=0$.)
For general points $p=(x,\ell),p'=(x',\ell')\in\Delta$,
define $s(p,q)=s(x,x')$ if $\ell=\ell'$ and $s(p,q)=0$ otherwise.
This defines a separation time $s:\Delta\times\Delta\to\N$.
We have the projection $\pi:\Delta\to M$ given by $\pi(x,\ell)=T^\ell x$
and satisfying $\pi T=F\pi$.

\begin{itemize}
\item[(A2)]   There is a distinguished subset or ``unstable leaf''
$W^u\subset\Lambda$ such that
each stable disk intersects $W^u$ in precisely one point, and 
there exist constants $C\ge1$, $\alpha\in(0,1)$ such that 
\begin{itemize}
\item[(i)] $d(T^nx,T^ny)\le C\alpha^n$, for all $y\in W^s(x)$, all $n\ge0$, and
\item[(ii)] $d(T^nx,T^ny)\le C\alpha^{s(x,y)}$ for all $x,y\in W^u$
and all $0\le n<R$.
\end{itemize}
\end{itemize}

\begin{rmk}  We note that Young~\cite{Young98} uses a separation time $s_0$
defined in terms of the underlying diffeomorphism $T:M\to M$ whereas our
separation time $s$ is defined in terms of the induced map 
$f:\Lambda\to\Lambda$.   
In particular,~\cite[conditions~(iii) and~(iv), p.~589]{Young98}
guarantee that $s_0\ge s$ and moreover that $s_0-(R-1)\ge s$.
Hence~\cite[assumption~(P4)(a)]{Young98} 
($d(T^nx,T^ny)\le C\alpha^{s_0(x,y)-n}$ for $0\le n<s_0(x,y)$) 
implies our assumption~(A2)(ii).

There is also a separation time in~\cite{Young98} that is denoted $s$.
This is different from our separation time and plays no role in this paper.
\end{rmk}

Let $\bar\Lambda=\Lambda/\sim$ where $x\sim x'$ if $x\in W^s(x')$.
Similarly, define the partition $\{\bar\Lambda_j\}$ of $\bar\Lambda$.
We obtain a well-defined return time function  $R:\bar\Lambda\to\Z^+$ and
induced map $f:\bar\Lambda\to\bar\Lambda$.
Let $F:\bar\Delta\to\bar\Delta$ denote the corresponding suspension map.
We note that this can be viewed as the quotient of $F:\Delta\to\Delta$
where $(x,\ell)$ is identified with $(x',\ell')$ if $\ell=\ell'$
and $x'\in W^s(x)$.
Let $\bar\pi:\Delta\to\bar\Delta$ denote the natural projection.

The separation time on $\Delta$ drops down to a separation time on $\bar\Delta$
(and agrees with the natural separation time defined using 
$f:\bar\Lambda\to\bar\Lambda$ and the partition $\{\bar\Lambda_j\}$).

\begin{itemize}
\item[(A3)]  The map $f:\bar\Lambda\to\bar\Lambda$ and partition $\{\bar\Lambda_j\}$ separate points in $\bar\Lambda$.
\end{itemize}
It follows that $d_\theta(p,q)=\theta^{s(p,q)}$ defines 
a metric on $\bar\Delta$ for each $\theta\in(0,1)$.

\begin{itemize}
\item[(A4)]  There exist
$F$-invariant probability measures $\tilde m$ on $\Delta$ and $\bar m$ on 
$\bar\Delta$ such that
\begin{itemize}
\item[(i)]  $\pi:\Delta\to M$ and $\bar\pi:\Delta\to \bar\Delta$ are
measure-preserving ($\pi$ takes $\tilde m$ to $m$ and $\bar\pi$ takes
$\tilde m$ to $\bar m$); and
\item[(ii)]  $F:\bar\Delta\to\bar\Delta$ is a Young tower (in the sense 
of section~\ref{sec-expand}(d)).
\end{itemize}
\end{itemize}

We say that an observation $\psi:\Delta\to\R$ {\em depends only on
future coordinates} if $\psi(p)=\psi(q)$ whenever $p\sim q$
where $\sim$ is the equivalence relation on $\Delta$ arising from quotienting
along stable disks.  Such an observation drops down to an
observation $\psi:\bar\Delta\to\R$.
The following result shows that any H\"older observation on $M$
is related to a Lipschitz observation on $\bar\Delta$
(cf.\ \cite{Sinai72,Bowen75}).

\begin{lemma} \label{lem-future}
Suppose that $\phi:M\to\R$ is $\gamma$-H\"older with respect to the metric $d$.
Then there exist functions $\psi,\chi:\Delta\to\R$ such that
\begin{itemize}
\item[(i)] $\phi\circ\pi=\psi+\chi-\chi\circ F$,
\item[(ii)] $\chi$ is bounded,
\item[(iii)] $\psi$ depends only on future coordinates,
\item[(iv)] $\psi:\bar\Delta\to\R$ is Lipschitz with respect to the
metric $d_\theta$, for $\theta=\alpha^{\gamma/2}$.
\end{itemize}
\end{lemma}

\begin{proof}
Given $p=(x,\ell)\in\Delta$, define $\hat p=(\hat x,\ell)$ where
$\hat x$ is the unique point in $W^s(x)\cap W^u$ (see (A2)).
Define 
\[
\chi (p)=\sum_{j=0}^{\infty}\phi (\pi F^jp)-\phi(\pi F^j \hat p).
\]
Note that $\pi F^jp=T^j\pi p=T^{j+\ell}x$
and similarly $\pi F^j\hat p=T^{j+\ell}\hat x$.
Since $x$ and $\hat x$ lie in the 
same stable disk $W^s$, it follows from (A2)(i) that
\begin{align*}
|\chi(p)| & 
\le \sum_{j=0}^\infty |\phi(\pi F^jp)-\phi(\pi F^j\hat p) 
\le  |\phi|_\gamma \sum_{j=0}^\infty d(T^{j+\ell}x,T^{j+\ell} \hat x)^\gamma \\
& \le |\phi|_\gamma C^\gamma \sum_{j=0}^\infty 
\alpha^{j\gamma} = |\phi|_\gamma C^\gamma (1-\alpha^\gamma)^{-1}.
\end{align*}

Define $\psi=\phi\circ\pi -\chi + \chi\circ F$.  Then
$\psi(p)=\sum_{j=0}^{\infty} \phi (\pi F^j\hat p)-\phi (\pi F^j\hat{ Fp})$
depends only upon future coordinates. 
It remains to check that $\psi$ is Lipschitz with respect to the metric 
$d_{\theta}$.   In fact, we prove that $\psi$ is Lipschitz with
respect to $d_{\theta^{1/2}}$ where $\theta=\alpha^\gamma$.

For any $N\ge1$, $p,q\in\Delta$,
\begin{align} \label{eq-chi} 
& |\psi (p)-\psi (q)|\le \sum_{j=0}^N |\phi (\pi F^j \hat p)-\phi (\pi F^j \hat q)|
+\sum_{j=0}^{N-1}|\phi(\pi F^j \hat{Fp})-\phi (\pi F^j \hat{Fq})|\\
&\qquad \qquad + \sum_{j=N+1}^{\infty} |\phi (\pi F^j \hat p)-\phi(\pi F^{j-1} \hat{Fp})|
+ \sum_{j=N+1}^{\infty} |\phi (\pi F^j \hat q)-\phi(\pi F^{j-1} \hat{Fq})|.
\nonumber
\end{align}
Suppose that $d_\theta(p,q)=d_\theta(\hat p,\hat q)\approx\theta^{2N}$.   
We show that each of these four terms is bounded by 
$\theta^N\approx  d_{\theta^{1/2}}(p,q)$ up to a constant.

Starting with the third term in~\eqref{eq-chi}, we note that 
$F\hat p=\hat{Fp}$ unless $p=(x,R(x)-1)$,
in which case $F\hat p=(f\hat x,0)$ and $\hat{Fp}=(\hat{fx},0)$.
Then $\pi F^j\hat p=T^{j-1}(f\hat x)$ and 
$\pi F^{j-1}\hat{Fp}= T^{j-1}(\hat{fx})$.
Since $f\hat x$ and $\hat{fx}$ lie in the same stable disk $W^s$, we have 
$|\phi (\pi F^j \hat p)-\phi (\pi F^{j-1} \hat{Fp})| 
\le |\phi|_\gamma C^\gamma \alpha^{(j-1)\gamma}$
so that $\sum_{j=N+1}^\infty
|\phi (\pi F^j \hat p)-\phi (\pi F^{j-1} \hat{Fp})| \le C'\theta^N$
as required.  Similarly for the fourth term in~\eqref{eq-chi}.  

Next, we consider the first term in~\eqref{eq-chi}.
By assumption, $s(p,q)\approx 2N$ so separation does not takes place during
the calculation.
Write $p=(x,\ell)$, $q=(y,\ell)$.   Then 
$\pi F^j\hat p=T^{j+\ell}\hat x=T^L f^J \hat x$
where $J\le j$ and $L<R(f^J\hat x)$.
Similarly, $\pi F^j\hat q=T^L f^J \hat y$.
Hence by (A2)(ii),
\begin{align*}
|\phi(\pi F^j \hat p)-\phi(\pi F^j \hat q)| & \le |\phi|_\gamma
d\bigl(T^L f^J \hat x,T^L f^J \hat y\bigr)^\gamma 
 \le |\phi|_\gamma C^\gamma \alpha^{s(f^J \hat x,f^J \hat y)\gamma}  \\ &
 = |\phi|_\gamma C^\gamma \alpha^{[s(\hat x,\hat y)-J]\gamma}  
 \le  |\phi|_\gamma C^\gamma \alpha^{[s(\hat x,\hat y)-j]\gamma}  
\approx |\phi|_\gamma C^\gamma\theta^{2N-j},
\end{align*}
so that 
$\sum_{j=0}^N |\phi (\pi F^j \hat p)-\phi (\pi F^j \hat q)| 
\le C'\theta^N$ 
as required.  Similarly for the second term in~\eqref{eq-chi}.  
\end{proof}

\begin{rmk}   Although Lemma~\ref{lem-future} is modelled on the treatments
in~\cite{Bowen75,ParryPoll90}, we have not defined a metric on $\Delta$
and hence the usual regularity statement about $\chi$ is missing.
\end{rmk}

\begin{thm} \label{thm-hyp}  
 Suppose that $T:M\to M$ satisfies (A1)--(A4) and
assume that $R\in L^{2+\delta}(\Lambda)$ for some $\delta>0$.
Let $\phi:M\to\R$ be a mean zero H\"older observation.
Then $\phi$ satisfies the ASIP.
\end{thm}

\begin{proof}   
Since $\pi:\Delta\to M$ is measure preserving, it suffices to prove the
ASIP for the lift $\tilde\phi=\phi\circ\pi:\Delta\to\R$.
By Lemma~\ref{lem-future}, there exists $\psi:\Delta\to\R$
depending only on future coordinates such that
$\tilde\phi_N-\psi_N$ is uniformly bounded, and it suffices to prove
the ASIP for $\psi$.   Since the projection $\bar\pi:\Delta\to\bar\Delta$
is measure preserving, it suffices to prove the ASIP for $\psi$
at the level of $\bar\Delta$.
Finally, Lemma~\ref{lem-future} guarantees that $\psi:\bar\Delta\to\R$
 is Lipschitz with respect to $d_\theta$,
so it suffices to prove the ASIP for Lipschitz observations on $\bar\Delta$
which is a Young tower by (A4)(ii).
Now apply Theorem~\ref{thm-ASIPtower}.~
\end{proof}

\paragraph{Nonuniformly hyperbolic flows}
Given an $L^1$ roof function $h:M\to\R^+$, we define a suspension flow
$T_t:M^h\to M^h$ in the same way that we defined the semiflow
in Section~\ref{sec-expand}(e).   If $T:M\to M$ satisfies (A1)--(A4),
we say that $T_t:M^h\to M^h$ is a {\em nonuniformly hyperbolic flow}.

\begin{cor} \label{cor-flow} Let $T_t:M^h\to M^h$ be a nonuniformly
hyperbolic flow.   Assume moreover that 
the return time function $R$ lies in $L^{2+\delta}(Y)$ and that
the roof function $h:M\to\R^+$ is H\"older.
Let $\psi:M^h\to\R$ be a mean zero H\"older observation.
Then $\psi$ satisfies the ASIP.
\end{cor}

\begin{proof}
This follows immediately from Theorem~\ref{thm-hyp}, applying Theorem~\ref{thm-MT}.
\end{proof}

\begin{rmk} \label{rmk-flow}
The weakened hypotheses mentioned in Remark~\ref{rmk-semiflow} apply
equally in the nonuniformly hyperbolic setting.
\end{rmk}

\section{Applications}
\label{sec-app}

In this section, we indicate a wide range of applications to which
the results in this paper apply.

We begin with nonuniformly expanding systems that can be modelled
by a Young tower as in Section~\ref{sec-expand}.
In the literature it is standard to speak of {\em return time asymptotics}
in the form $m\{y\in Y:R(y)\ge n\}=O(n^{-\gamma})$.
(Recall from Section~\ref{sec-expand} that $Y$ is the subset used for inducing,
equivalently the base of the Young tower.)

\begin{prop} \label{prop-R}
If $m\{R\ge n\}=O(n^{-\gamma})$ for some $\gamma>2$,
then $R\in L^{2+\delta}(Y)$ for $\delta\in(0,\gamma-2)$.
\end{prop}

\begin{proof}
This is immediate from the inequality
$E[R^{2+\delta}]\le \sum_{n=0}^{\infty}m\{ R^{2+\delta}\ge n\}
= \sum_{n=0}^{\infty}m\{ R\ge n^{\frac{1}{2+\delta}}\}$.
\end{proof}

Many maps satisfy the condition in Proposition~\ref{prop-R}:

\vspace{1ex}
\noindent (i) the Alves-Viana map~\cite{ALP02} $T:S^1\times I\to S^1\times I$
\begin{align*}
T(\omega,x)=(16\omega,a-x^2+\epsilon\sin(2\pi \omega))
\end{align*}
when $0$ is preperiodic for the map $x\mapsto a-x^2$ and
$\epsilon$ is small enough.

\vspace{1ex}
\noindent (ii) the Liverani-Saussol-Vaienti (Pomeau-Manneville) 
maps~\cite{LSV99} $T:[0,1]\to[0,1]$
\[
Tx=\left\{ \begin{array}{cc}
x(1+2^{\alpha}x^{\alpha})& 0\le x<\frac{1}{2}\\
2x-1 &\frac{1}{2}\le x<1
\end{array}\right.
\]
for $0<\alpha<\frac{1}{2}$.

\vspace{1ex}
\noindent (iii)  certain classes of multimodal maps, 
Bruin~{\em et al.}~\cite{BLS}.

\vspace{1ex}
\noindent (iv)  
 a class of expanding circle maps $T:S^1\to S^1$ of degree $d>1$
with a neutral fixed point, Young~\cite[Section 6]{Young99}:
$T$ is $C^1$ on $S^1$ and $C^2$ on $S^1-\{0\}$, $T'>1$ on $S^1-\{0\}$,
$T(0)=0$, $T'(0)=1$, and for $x\not =0$,
$-xT''(x)\simeq |x|^{\alpha}$ for $0<\alpha <\frac{1}{2}$.

\vspace{1ex}
\noindent  Applying Theorem~\ref{thm-intro}, we obtain the ASIP for
H\"older observations for the systems in (i)--(iv) above.
For example, in (iii) and (iv) we obtain the ASIP under the same conditions 
for which~\cite{BLS} and~\cite{Young99} obtain the CLT.

\vspace{1ex}

Next, we recall examples of nonuniformly hyperbolic systems that have been
modelled by towers.
Consider the following classes  of $C^{1+\epsilon}$ diffeomorphisms 
treated in Young~\cite{Young98} (see also Baladi~\cite[\S4.3]{Baladi}):

\vspace{1ex}
\noindent(v) Lozi maps and certain  piecewise hyperbolic maps 
\cite{Young98,Chernov99a}.

\vspace{1ex}
\noindent(vi) a class of H\'{e}non 
maps~\cite{Benedicks-Young1,Benedicks-Young2}.

\vspace{1ex}
\noindent(vii) some partially hyperbolic diffeomorphisms with a mostly
contracting direction~\cite{Castro,Dolgopyat}.

\vspace{1ex}
In these examples, the return time asymptotics are exponential so
certainly $R\in L^{2+\delta}$.
By Theorem~\ref{thm-hyp}, we obtain the  
ASIP for H\"{o}lder observations for the systems in (v)--(viii) above.

\paragraph{Billiard maps and Lorentz flows}
Finally, we consider the application to the planar periodic
Lorentz gas discussed in the introduction.
Under the finite horizon condition,
Young~\cite{Young98} demonstrated that the billiard map
(which is the Poincar\'e map for the flow) is nonuniformly
hyperbolic with exponential return time asymptotics.
As a result, Young established exponential decay of correlations
for such billiard maps, resolving a long-standing (and controversial)
open question.   Chernov~\cite{Chernov99b} extended Young's method
to obtain the same result for infinite horizons.

For our purposes, the weaker conclusion that $R\in L^{2+\delta}$ is again 
sufficient.  Hence, by the results in~\cite{Chernov99b,Young98},
the first statement of Theorem~\ref{thm-Lorentz}
is an immediate consequence of Theorem~\ref{thm-hyp}.

For the flow itself, the finite horizon condition is crucial since even the
CLT is unlikely in the infinite horizon case.
Assuming finite horizons, the roof function $h$ is uniformly bounded
and piecewise H\"older.  Since $h$ is not uniformly H\"older,
Corollary~\ref{cor-flow} does not apply directly, but the result is
easily modified as in Remarks~\ref{rmk-semiflow} and~\ref{rmk-flow}
to include such roof functions.
Hence, we obtain the second statement of Theorem~\ref{thm-Lorentz}.


\appendix
\section{ASIP for functions of mixing sequences}

Here is a special case of Philipp \& Stout~\cite[Theorem 7.1]{PhilippStout75} 
adapted to dynamical systems terminology.
The notation is as in Section~\ref{sec-expand}(c).

\begin{thm}[Philipp \& Stout] \label{thm-PS}
Assume that there exists $\delta\in(0,2]$, $\sigma^2>0$ and $C>0$ such that
for all $k,N\ge1$,
\begin{itemize}
\item[(i)] $v\in L^{2+\delta}(\Lambda)$ and $\int_\Lambda v\,dm=0$,
\item[(ii)] $\int_\Lambda v_N^2 \,dm = \sigma^2N+O(N^{1-\delta/30})$,
\item[(iii)] $\sum_{a\in\alpha_0^{k-1}}\int_a|v-\frac{1}{m(a)}\int_a v\,dm|^{2+\delta}dm
\le Ck^{-(2+7/\delta)(2+\delta)}$,
\item[(iv)] $\bigl|m(a\cap f^{-(N+k)}(b))-m(a)m(b)\bigr|\le CN^{-168(1+2/\delta)}$
for all $a\in\alpha_0^{k-1}$ and all measurable sets $b$.
\end{itemize}
Then $v_N=W(N)+O(N^{1/2-\delta/600})$.
\end{thm}

\section{ASIP for suspensions}

Suppose that $(\Lambda,m)$ is a probability space
and that $f:\Lambda\to\Lambda$ is a measure-preserving transformation.
Let $h:\Lambda\to\R^+$ be a roof function and suppose that
$f_t:\Lambda^h\to\Lambda^h$ is the corresponding suspension (semi)flow
as in Section~\ref{sec-expand}(e).   The following result is a special case
of~\cite[Theorem~4.2]{MT02}.

\begin{thm}[Melbourne \& T\"or\"ok] \label{thm-MT}
Let $\delta>0$.
Suppose that $h\in L^{2+\delta}(\Lambda)$ and that
$\sum_{j=0}^{N-1}h\circ f^j=N\bar h+o(N^{1-\delta})$ as $N\to\infty$
almost everywhere.   

Suppose that
$\psi:\Lambda^h\to\R$ lies in $L^\infty(\Lambda^h)$ and has mean zero.
Define $\phi:\Lambda\to\R$ by $\phi(x)=\int_0^{h(x)}\psi(f_tx)$.
If $\phi$ satisfies the ASIP on $\Lambda$ with variance $\sigma_1^2$,
then $\psi$ satisfies the ASIP on $\Lambda^h$ with variance 
$\sigma^2=\sigma_1^2/\bar h$.
\end{thm}

Theorem~\ref{thm-MT} is easily modified for discrete suspensions.
Let $R:\Lambda\to\Z^+$ be an $L^1$ return time function and form
the discrete suspension map
$F:\Delta\to\Delta$ as in Section~\ref{sec-expand}(d).

\begin{cor} \label{cor-MT}
Let $\delta>0$.
Suppose that $R\in L^{2+\delta}(\Lambda)$ and that
$\sum_{j=0}^{N-1}R\circ f^j=N\bar R+o(N^{1-\delta})$ as $N\to\infty$
almost everywhere.   

Suppose that
$\phi:\Delta\to\R$ lies in $L^\infty(\Delta)$ and has mean zero.
Define $\Phi:\Lambda\to\R$ by $\Phi(x)=\sum_{j=0}^{R(x)-1}\phi(f^jx)$.
If $\Phi$ satisfies the ASIP on $\Lambda$ with variance $\sigma_1^2$,
then $\phi$ satisfies the ASIP on $\Delta$ with variance 
$\sigma^2=\sigma_1^2/\bar R$.
\end{cor}

 \paragraph{Acknowledgements}  
This research was supported in part by EPSRC Grant GR/S11862/01.
IM is greatly indebted to UH for the use of e-mail, given that pine is 
currently not supported on the University of Surrey network.

\end{document}